\documentclass{IEEEtran4PSCC}
\ifCLASSINFOpdf
   \usepackage[pdftex]{graphicx}
\else
   \usepackage[dvips]{graphicx}
\fi
\usepackage[cmex10]{amsmath}
\usepackage{amssymb,amsthm,epsfig,dsfont,color,subfigure,empheq,graphicx,balance,multirow}
\usepackage{enumerate,url,algpseudocode,algorithm,wasysym,soul,bbm}
\usepackage[table]{xcolor}
\renewcommand{\d}[1]{\ensuremath{\operatorname{d}\!{#1}}}

\newcommand \bzero{\mathbf{0}}

\newcommand \be{\mathbf{e}}
\newcommand \bg{\mathbf{g}}

\newcommand \bdell{\boldsymbol{\ell}} 

\newcommand \bu{\mathbf{u}}

\newcommand \bw{\mathbf{w}}
\newcommand \bx{\mathbf{x}}
\newcommand \by{\mathbf{y}}
\newcommand \bz{\mathbf{z}}

\newcommand \bB{\mathbf{B}}

\newcommand \bG{\mathbf{G}}

\newcommand \bI{\mathbf{I}}

\newcommand \bY{\mathbf{Y}}

\newcommand \blambda{\boldsymbol{\lambda}}

\newcommand \bpi{\boldsymbol{\pi}}

\newcommand \bphi{\boldsymbol{\phi}}

\newcommand \mcE{\mathcal{E}}

\newcommand \mcN{\mathcal{N}}

\newcommand \mcX{\mathcal{X}}

\newcommand \bby{\bar{\mathbf{y}}}

\makeatletter
\let\old@ps@headings\ps@headings
\let\old@ps@IEEEtitlepagestyle\ps@IEEEtitlepagestyle
\def\psccfooter#1{%
    \def\ps@headings{%
        \old@ps@headings%
        \def\@oddfoot{\strut\hfill#1\hfill\strut}%
        \def\@evenfoot{\strut\hfill#1\hfill\strut}%
    }%
    \def\ps@IEEEtitlepagestyle{%
        \old@ps@IEEEtitlepagestyle%
        \def\@oddfoot{\strut\hfill#1\hfill\strut}%
        \def\@evenfoot{\strut\hfill#1\hfill\strut}%
    }%
    \ps@headings%
}

\allowdisplaybreaks

\begin{document}
%
\title{DNN-based Policies for Stochastic AC OPF}

\author{
    \IEEEauthorblockN{Sarthak Gupta\IEEEauthorrefmark{1}\IEEEauthorrefmark{2}, Sidhant Misra\IEEEauthorrefmark{2}, Deepjyoti Deka\IEEEauthorrefmark{2}, Vassilis Kekatos\IEEEauthorrefmark{1}}
    \IEEEauthorblockA{\IEEEauthorrefmark{1}ECE Dept., Virginia Tech, Blacksburg,  VA 24061, USA\\
    \{gsarthak,kekatos\}@vt.edu}
    \IEEEauthorblockA{\IEEEauthorrefmark{2}Theory Division, 
 Los Alamos National Laboratory, Los Alamos, NM, USA\\
  \{sidhant,deepjyoti\}@lanl.gov}
}
\mark{Gupta, Kekatos, and Jin: Controlling Smart Inverters using Proxies: A Chance-Constrained DNN-based Approach}

\maketitle

\begin{abstract}
A prominent challenge to the safe and optimal operation of the modern power grid arises due to growing uncertainties in loads and renewables. Stochastic optimal power flow (SOPF) formulations provide a mechanism to handle these uncertainties by computing dispatch decisions and control policies that maintain feasibility under uncertainty. Most SOPF formulations consider simple control policies such as affine policies that are mathematically simple and resemble many policies used in current practice. 
Motivated by the efficacy of machine learning (ML) algorithms and the potential benefits of general control policies for cost and constraint enforcement, we put forth a deep neural network (DNN)-based policy that predicts the generator dispatch decisions in real time in response to uncertainty. The weights of the DNN are learnt using stochastic primal-dual updates that solve the SOPF without the need for prior generation of training labels and can explicitly account for the feasibility constraints in the SOPF. The advantages of the DNN policy over simpler policies and their efficacy in enforcing safety limits and producing near optimal solutions are demonstrated in the context of a chance constrained formulation on a number of test cases. 
\end{abstract}


\section{Introduction}
In current intra-day power systems operations, the optimal power flow (OPF) is solved at a time scale of 5-15 minutes using load forecasts to obtain economic generation dispatch, whereas affine/linear generation control policies are used within each OPF time period to account for real-time fluctuations in generation and demand. With growing levels of uncertainty, the design of the real-time control policies have become increasingly important for secure and economic grid operation. Such control design is the primary focus of this paper. 

Motivated by the need for better uncertainty management, the stochastic AC-OPF (SOPF) \cite{bienstock2014chance,roald2013analytical} problem that accounts for uncertainty, has received significant attention from the research community. Being a stochastic extension of the non-convex AC-OPF, the stochastic AC-OPF is highly computationally challenging and much effort has been devoted towards obtaining tractable formulations and designing efficient algorithms for obtaining fast and reliable solutions to this problem. The design of control polices however is relatively under-studied and most SOPF formulations consider \emph{affine} functions to represent control policies. These linearly parameterized policies are algorithmically easier to handle, and serve as reasonable mathematical models for the automatic generation control and local voltage control employed in current practice. However, affine control polices are restrictive, and more general control policies need to be studied as they have the potential to  handle a much larger class of real-time fluctuations, possibly with limited information.  

Design of general control policies poses significant technical/computational challenges and the choice of parameterization for their mathematical representation plays a pivotal role in their design. First, the expressive power of the parameterization dictates how general of a control policy it can represent. Second, the parameterization chosen must be compatible with the corresponding SOPF formulation to ensure computational tractability. In what follows, we provide a broad classification of approaches in the literature to this coupled control-parameterization and algorithm-design problem:\\
\emph{i) Affine policy.} Here the control is restricted to be a linear function typically of the net system uncertainty, which makes it easy to model within OPF. Such affine policies have been used for convex SOPF with the DC power flow model \cite{bienstock2014chance} as well as for scenario-based approaches \cite{vrakopoulou2013probabilistic}, and SOPF with nonlinear AC power flow model \cite{dall2017chance,MEZGHANI2020}. \\
\emph{ii) Non-affine policies} include non-linear and hence more general policies that can ensure better feasibility and optimality enforcement over the affine ones. However they are computationally harder to incorporate within the SOPF. Examples include polynomial chaos based policies \cite{muhlpfordt2019chance,metivier2020efficient}, as well as non-liner policies with generator saturation \cite{roald2015optimal}.\\
\emph{iii) Nonlinear policies using robust optimization.} This approach uses a fully optimized recourse for each uncertainty realization using a min-max-min formulation \cite{lorca2017robust}. The problem is highly computationally challenging and approximation and relaxations are employed that can lead to conservative solutions. Similar robust formulations with affine control have also been investigated in the literature \cite{lee2021robust}.\\
\emph{iv) Nonlinear policies using data-driven methods.} This approach uses historical data or simulated OPF solutions to devise efficient non-linear control policies as well as to improve the computational speed of the SOPF. Examples of such methods include kernel-based control policies for voltage control \cite{JKGD19} and active set learning approaches \cite{ng2018statistical,DekaMisraPowerTech19} for OPF. Owing to the advancements in deep learning over the past decade, DNN-based OPF solution schemes have also been explored by the power systems community. Taking the conventional supervised learning route, DNN-based OPF solvers have been trained to match the optimal labels~\cite{zamzam2020learning,chen2020learning}. However, solving a large number of OPFs to generate optimal labels is itself a computationally overwhelming task. Sensitivity-informed deep learning~\cite{SGKCB2020,L2OPWRS21}, which matches not only OPF minimizers but also partial derivatives with respect to the inputs, partially alleviates the above concern by improving upon the data efficiency. Nonetheless, none of these DNN solvers provides sound mechanisms to incorporate constraints or handle randomness.

In this paper, we consider a  DNN-based control policy design for SOPF. The DNN parameterization serves as a generalization to non-linear policies, but also provides backward compatibility with previous policies since it can easily accommodate communication constraints and restricted features (e.g., control based on total/net uncertainty and/or local control). Different from the aforementioned works, our approach does not involve an offline data generation phase where a large number of OPF or OPF-like problems are solved. Instead, similar to~\cite{SG_SGC20,OPFandLearnTSG21}, the training phase of the DNN itself serves as the SOPF solver, which if deployed within current practices will replace the OPF solved at the 5-15 minutes time scale. Crucially, we also model the non-linear power system constraints (both hard and stochastic) during training to ensure feasibility of the control actions. Thus our approach uses a hybrid model that has the advantages of the NNs such as high representation power and the potential for parallelization and GPU implementation, without sacrificing feasibility like black-box data-driven efforts.

The rest of the paper is organized as follows. Section~\ref{sec:problem} presents the mathematical formulation for SOPF with system and control constraints. Section~\ref{sec:DNNtrain} describes in detail our control policy design, in particular the modeling and training of our neural network with system constraints. Section~\ref{sec:numerical} includes experimental validation of our approach and comparison with existing methods on IEEE test cases. Finally, Section~\ref{sec:conclusion} concludes the paper.

\section{Problem Formulation}\label{sec:problem}
\subsection{Grid Modeling}\label{subsec:gridmodel}
Consider a power transmission network modeled as a directed connected graph $\mathcal{G}=(\mcN,\mathcal{E})$. The nodes $n \in \mcN:=\{1,\dots,N\}$ of the graph $\mathcal{G}$ correspond to network buses, while the directed edges $e_{nk} \in \mathcal{E}$ to transmission lines.
The complex voltage at bus $n$ is expressed in polar coordinates as $v_ne^{j\theta_n}$, and the complex power injection at the same bus as $p_n+jq_n$. The power injections at node $n$ are described by the AC power flow equations (AC-PF) as
\vspace{-.8em}
\begin{subequations}\label{eq:pf_inj}
\begin{align}
    p_n&=v_n\sum^{N}_{k=1}v_k\left(G_{nk}\cos\theta_{nk}+B_{nk}\sin\theta_{nk}\right)\\
    q_n&=v_n\sum^{N}_{k=1}v_n\left(G_{nk}\sin\theta_{nk}-B_{nk}\cos\theta_{nk}\right)
\end{align}
\end{subequations}
where $\theta_{nk}:=\theta_n-\theta_k$, and $G_{nk}$ and $B_{nk}$ are the entries of the real and imaginary parts of the bus admittance matrix $\bY=\bG+j\bB$. Power injections can be decomposed into their dispatchable $(p^g_n,q^g_n)$ and inflexible $(p^d_n,q^d_n)$ parts as 
\[p_n=p^g_n-p^d_n \quad \textrm{and} \quad q_n=q^g_n-q^d_n.\]
The former corresponds to generators and flexible loads; the latter to inelastic loads hosted per bus $n$. The buses hosting dispatchable injections form set $\mcN_ g\subset \mcN$, and its complement is defined as $\mcN_ \ell$. Set $\mcN_ g$ will be henceforth referred to as the set of \emph{generator} buses, and $\mcN_\ell$ as the set of \emph{load} buses. To simplify the exposition, each generator bus is assumed to host exactly one dispatchable unit. The bus indexed by $n=1$ is designated as the slack and reference bus with $\theta_1=0$. 

The complex power $P_{nk}+jQ_{nk}$ flowing from bus $n$ to bus $k$ over line $e_{nk}$ can be expressed as~
\begin{subequations}\label{eq:pf_flow}
\begin{align}%
    P_{nk}&=v_nv_k\left(G_{nk}\cos\theta_{nk}+B_{nk}\sin\theta_{nk}\right)-v_n^2G_{nk}\\
    Q_{nk}&=v_nv_k\left(G_{nk}\sin\theta_{nk}-B_{nk}\cos\theta_{nk}\right)+v_n^2\left(B_{nk}-B_{nk}^\text{sh}\right)
\end{align}
\end{subequations}
where $B_{nk}^\text{sh}$ is the shunt susceptance of line $e_{nk}$. The apparent power flow on the line is
\begin{equation}\label{eq:flow_VA}
f_{n,k}:=\sqrt{P_{n,k}^2+Q_{n,k}^2}.
\end{equation}


\subsection{Optimal Power Flow (OPF)}\label{subsec:OPF}


Let $c_n(p^g_n)$ denote the cost of generation at bus $n$. Given inflexible loads $\{p^d_n,q^d_n\}_{n\in\mcN}$, the OPF problem aims at minimizing the total cost of generation while operating the transmission network within limits. The OPF is posed as
\begin{subequations}\label{eq:opf}
\begin{align}
    \min~&~\sum_{n\in \mcN_g}c_n(p^g_n) \\
    \mathrm{over}~&\{v_n,\theta_n\}_{n \in \mcN},\{p^g_n,q^g_n\}_{n \in \mcN_g} \nonumber\\
    \mathrm{s.to}~&~\eqref{eq:pf_inj},\eqref{eq:pf_flow}&\forall n \in \mcN, \forall e_{nk} \in \mcE \nonumber\\
    &~p_n=p^g_n-p^d_n&\forall n \in \mcN_ g\label{eq:png}\\
    &~q_n=q^g_n-q^d_n&\forall n \in \mcN_ g\label{eq:qng}\\
    &~p_n=-p^d_n&\forall n \in \mcN_\ell\label{eq:pnl}\\
    &~q_n=-q^d_n&\forall n \in \mcN_\ell\label{eq:qnl}\\
    &~\underline{p}^g_n \leq p^g_n\leq \bar{p}^g_n &\forall n \in \mcN_g \label{eq:pglim}\\
    &~\underline{q}^g_n \leq q^g_n\leq \bar{q}^g_n &\forall n \in \mcN_g\label{eq:qglim}\\
    &~\underline{v}_n \leq v_n\leq \bar{v}_n &\forall n \in \mcN\label{eq:vlim}\\
    &~f_{n,k}\leq \bar{f}_{n,k} &\forall e_{n,k} \in \mathcal{E}\label{eq:flowlim}\\
    &~\theta_1=0\label{eq:refangle}
    \end{align}
\end{subequations}
where \eqref{eq:pf_inj} and \eqref{eq:png}--\eqref{eq:pnl} enforce the power flow equations for power injections. Constraints \eqref{eq:pglim}-\eqref{eq:qglim} represent generation limits. Constraint \eqref{eq:vlim} confines voltage magnitudes within given limits. Constraints \eqref{eq:pf_flow} and \eqref{eq:flowlim} ensure that apparent power flows remain within line ratings. Lastly, constraint \eqref{eq:refangle} fixes the phase angle at the reference bus.

Given load demands stored in vector $\bphi:=\{p^d_n,q^d_n\}_{n\in\mcN}$, the system operator solves \eqref{eq:opf} to find the optimal voltage and active power set-points for generators, which can be collected in vector $\bx:=\{\{v_n\}_{n\in\mcN_ g},\{p_n^g\}_{n\in\mcN_ g\setminus \{1\}}\}$. Given $(\bphi,\bx)$, all other grid quantities involved in the OPF can be expressed as functions of $(\bphi,\bx)$ implicitly via the AC-PF equations \eqref{eq:pf_inj} and \eqref{eq:pf_flow}. These quantities include the bus generations, and $(p^g_1,q_g^1)$ which we collectively represent using the variable $\by$. The OPF of \eqref{eq:opf} can be abstracted as a parametric optimization solely over variable $\bx$ given parameters $\bphi$:
\begin{subequations}\label{eq:opf2}
\begin{align}
    \min_{\bx}~&~\sum_{n\in \mcN_g}c_n\left(\bx,\bphi\right)\\
    \mathrm{s.to}~&~\underline{\bx}\leq\bx\leq\bar{\bx}\label{eq:xlim}\\
    &~ \by\left(\bx,\bphi\right)\leq\bar{\by}.\label{eq:ylim}
\end{align}
\end{subequations}
Here \eqref{eq:xlim} captures the constraints on generator setpoints, that is \eqref{eq:pglim} and \eqref{eq:vlim} for $n\in\mcN_g$, while \eqref{eq:ylim} captures the constraints in \eqref{eq:qglim}, \eqref{eq:vlim}, and \eqref{eq:flowlim}. 

Heed there is no need to explicitly enforce the constraints on loads as in \eqref{eq:png}--\eqref{eq:qnl}. This is because active and reactive load demands have been included in $\bphi$, and their effect on all other grid quantities of interest is captured by the power flow equations, which are still taken into account in~\eqref{eq:opf2} indirectly through the constraint function $\by(\bx,\bphi)$. Note also that even though the mapping $\by(\bx,\bphi)$ does not feature an analytical form, we will be able to compute its partial derivatives with respect to $\bx$ via the chain rule and the inverse function theorem later in Section~\ref{subsec:grad}. For now, let us explain how the OPF in \eqref{eq:opf2} can be modified to cope with uncertainty in $\bphi$.

\subsection{Chance-Constrained Optimal Power Flow}
Solving the non-convex problem of \eqref{eq:opf2} can be computationally taxing if $\bphi$ changes frequently. Also by the time \eqref{eq:opf2} is solved and the optimal setpoints are communicated to generators, demands $\bphi$ may have changed rendering $\bx$ obsolete. The above concerns can be accounted for by considering a stochastic version of the OPF~\cite{ergodic}, \cite{6965430}. We consider the chance constrained OPF (CC-OPF), that treats $\bphi$ as a random variable and seeks a control policy $ \bx(\bphi)$ that reacts to the realization of $\bphi$ to minimize the expected generation cost while satisfying network constraints with high probability
\begin{subequations}\label{eq:opf3}
\begin{align}
    \min_{\bx(\bphi) \in\mathcal{X}}~&~\mathbb{E}\left[\sum_{n\in \mcN_g}c_n\left(\bx,\bphi\right)\right] \\
    &~\mathrm{Pr}\left[y_i\left(\bx,\bphi\right)\leq\bar{y}_i\right]\geq 1-\alpha, &i=1:M \label{eq:chancehigh}
\end{align}
\end{subequations}
where the expectation $\mathbb{E}$ and the probability $\mathrm{Pr}$ are with respect to $\bphi$; and parameter $\alpha\in(0,1)$ controls the probability of constraint violation. The deterministic constraint \eqref{eq:xlim} has been abstracted as $\bx\in\mathcal{X}$. Using the indicator function, the chance constraints in \eqref{eq:chancehigh}  can be recast as 
\begin{align} \label{eq:opf4}
    \mathbb{E}\left[\mathbbm{1}\left(y_i\left(\bx,\bphi\right)-\bar{y}_i\right)\right]\geq 1-\alpha,  \ &i=1:M
\end{align}
where the indicator function $\mathbbm{1}(x)$ takes the values of $1$ and $0$ for $x\leq0$ and $x>0$, respectively. 
The problem \eqref{eq:opf3} is a variational optimization problem over the control policy $\bx(\bphi)$ and is generally intractable in its native form. Instead, one can parameterize the control policy as $\bx(\bphi) = \bpi\left(\bphi ; \bw \right)$, where $\bpi(.)$ is a chosen function of $\bphi$ determined by the parameters $\bw$.
Restricting the policy to take this parameterized form, \eqref{eq:opf3} can be written as the constrained stochastic minimization

\begin{subequations}\label{eq:opf5}
\begin{align}\label{eq:indhigh}
    \min_{\bw:\bpi_{\bw} \in\mathcal{X}}~&~\mathbb{E}\left[\sum_{n\in \mcN_g}c_n\left(\bpi_{\bw},\bphi\right)\right]\\
    \mathrm{s.to}~&~\mathbb{E}\left[\mathbbm{1}\left(y_i\left(\bpi_{\bw},\bphi\right)-\bar{y}_i\right)\right]\geq 1-\alpha, &i=1:M\nonumber
\end{align}
\end{subequations}
where the optimization is now over the policy parameters $\bw$, and the notation has been slightly abused by simplifying $\bpi\left(\bphi;\bw\right)$ to $\bpi_{\bw}$.  Given their advanced representation capabilities, DNNs constitute great candidates for modeling the policy $\bpi_{\bw}$. In the case of DNNs, the parameter vector $\bw$ collects the weighs and biases across all layers of the DNN. 
In general \eqref{eq:opf5} is an inner approximation of \eqref{eq:opf3} due to the restriction imposed on the policy. However, given the universal approximation properties of DNNs, \eqref{eq:opf5} is nearly identical to \eqref{eq:opf3} when the DNN depths is sufficient.

We next delineate how a DNN can be trained to solve the CC-OPF in \eqref{eq:opf5}. Unlike standard feed-forward DNNs, our approach involves an iterative primal-dual scheme that inherently satisfies the non-convex stochastic constraints involved in \eqref{eq:opf5}.

\section{Finding Optimal Policies}\label{sec:DNNtrain}
This section approximates the indicator function in \eqref{eq:opf5} with a logistic function; adopts a stochastic primal-dual scheme to learn a DNN to solve the CC-OPF. It explains how the required gradients can be computed, and provides an overview of the training and operational phase of the DNN.

\subsection{Approximating the Indicator Function}
The weights $\bw$ for a DNN-based policy can be learned by solving \eqref{eq:opf5}. However, computing the expectations in \eqref{eq:opf5} is challenging even if the pdf of $\bphi$ is known, given the non-linearities in the objective and constraints. Standard DNN training alleviates this issue by using stochastic gradient descent (SGD) updates over a training dataset generated by drawing random samples with respect to the distribution of $\bphi$ to learn the DNN weights. In the presence of constraints, stochastic primal-dual updates (SPD) can be used to yield a similar training process~\cite{Ribeiro19},~\cite{OPFandLearnTSG21}. Unfortunately, the indicator function $\mathbbm{1}(x)$ in \eqref{eq:opf5} prevents the application of any gradient-based approach. This is because $\mathbbm{1}(x)$ is discontinuous at $x=0$, and its derivative becomes $0$ at $x\neq 0$. Hence, no useful gradients are obtained while applying SPD. 

\begin{figure}[t]
    \centering
    \hspace*{-2em}\includegraphics[scale=0.45]{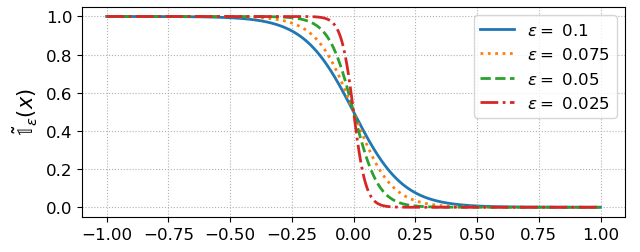}
    \hspace*{-2em}\includegraphics[scale=0.45]{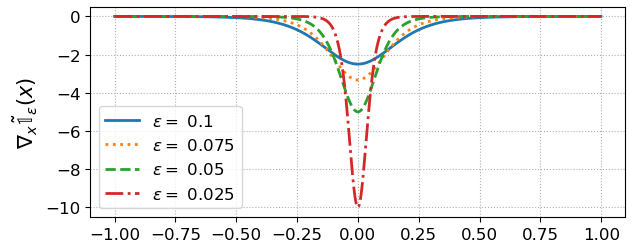}
    \caption{Logistic function approximation of the indicator function \emph{(top)}, and its derivative \emph{(bottom)} for different values of parameter $\epsilon$.}
    \label{fig:log}
    \vspace{-2em}
\end{figure}


The non-differentiability of $\mathbbm{1}(x)$ has been addressed in the literature by substituting it with a differentiable approximation. Depending on the objective and constraints, the convex approximations of $\mathbbm{1}(x)$ proposed in \cite{Nemirovski07} may yield overall convex CC-OPF formulations~\cite{dall2017chance}. Nonetheless, the computational advantage gained by introducing convexity is balanced off by the sub-optimal, conservative nature of these formulations. It is worth stressing that the problem in \eqref{eq:opf5} has sources of non-convexity other than $\mathbbm{1}(x)$, namely the policy $\bpi\left(\bphi;\bw\right)$ and the underlying non-linear power flow equations. Therefore, aiming for a more accurate and differentiable approximation of $\mathbbm{1}(x)$ is pertinent, even if this approximation is non-convex. Fortunately, the logistic function, which is well known in the ML community, can serve as a smooth surrogate of $\mathbbm{1}(x)$~\cite{Chen95}. The logistic function and its derivative are 
\begin{equation}
\tilde{\mathbbm{1}}_{\epsilon}(x):=\frac{e^{-x/\epsilon}}{1+e^{-{x}/{\epsilon}}},~~
\nabla_{x}\tilde{\mathbbm{1}}_{\epsilon}(x)=\frac{-\tilde{\mathbbm{1}}_{\epsilon}(x)\left(1-\tilde{\mathbbm{1}}_{\epsilon}(x)\right)}{\epsilon}\label{eq:gradlog}
\end{equation}
where parameter $\epsilon$ controls the accuracy of the approximation. It is easy to verify that $\tilde{\mathbbm{1}}_{\epsilon}(x)$ tends to $\mathbbm{1}_{\epsilon}(x)$ as $\epsilon\to 0^+$ as demonstrated in Figure~\ref{fig:log}. The same figure points towards a possible trade-off between approximation error and rate of convergence for SPD: As $\epsilon$ decreases, the range of values of $x$ over which $\tilde{\mathbbm{1}}_{\epsilon}(x)$ has non-diminishing derivative decreases as well. It is hence anticipated that for smaller values of $\epsilon$, a larger number of SPD updates may be needed.

\subsection{Stochastic Primal-Dual Updates}\label{subsec:SPD}
Dualizing \eqref{eq:opf5} with $\mathbbm{1}(x)$ being replaced by $\tilde{\mathbbm{1}}_{\epsilon}(x)$ provides the Lagrangian function
\begin{align}
    L(\bw;\blambda)&:=\mathbb{E}\left[\sum_{n\in \mcN_g}c_n\left(\bpi_{\bw},\bphi\right) \right]\nonumber\\
    &~~+\blambda^{\top}\left[(1-\alpha)\mathbf{1}-\mathbb{E}\left[\tilde{\mathbbm{1}}_{\epsilon}\left(\by\left(\bpi_{\bw},\bphi\right)-\bar{\by}\right)\right]\right]\label{eq:lagrang}
\end{align}
where $\blambda$ is the vector of Lagrange multipliers associated with the constraints in \eqref{eq:indhigh}. Observe that $\tilde{\mathbbm{1}}_{\epsilon}\left(\by\left(\bpi_{\bw},\bphi\right)-\bar{\by}\right)$ are vector quantities obtained by applying $\tilde{\mathbbm{1}}_{\epsilon}(x)$ element-wise for all constraints in \eqref{eq:indhigh}. The corresponding dual problem is
\begin{align}\label{eq:dual}
D^*:=\max_{\blambda\geq \bzero}\min_{\bw:\bpi_{\bw}\in\mcX}. L(\bw;\blambda)
\end{align}
A stationary point of this min/max problem can be reached via the projected primal-dual iterations indexed by $k$
\begin{subequations}\label{eq:pdupdate}
\begin{align}
    \bw_{k+1}&:=\big[\bw_{k}-\mu\nabla_{\bw}L(\bw_{k}\blambda_k)\big]_{\bpi_{\bw}\in\mcX}\label{eq:pdupdate:p}\\
    \blambda_{k+1}&:=\big[\blambda_{k}+\nu\nabla_{\blambda}L(\bw_{k};\blambda_{k})\big]_+\label{eq:pdupdate:dhigh}
\end{align}
\end{subequations}
with positive step sizes $\mu$ and $\nu$. The primal update \eqref{eq:pdupdate:p} involves a projection of $\bpi_{\bw}$ into the feasible set $\mcX$. With $\bpi_{\bw}$ being modeled by a DNN, projection onto $\mcX$ can be easily accomplished by using appropriately scaled hyperbolic tangent functions ($\tanh$). Dual variables are projected on the non-negative orthant via the operator $[x]_{+}:=\max\{x,0\}$.
 
The expectation operator in \eqref{eq:lagrang} can be handled using stochastic approximation. This entails first surrogating each of the expected values by their sample averages over $K$ uncertainty realizations $\{\bphi^k\}^K_{k=1}$. For the cost function for example, this approximation yields
\begin{align*}
  \mathbb{E}\left[\sum_{n\in \mcN_g}c_n\left(\bpi_{\bw},\bphi\right)\right]\simeq \frac{1}{K}\sum_{k=1}^{K}\sum_{n\in \mcN_g}c_n\left(\bpi_{\bw},\bphi^k\right).
\end{align*}
The \emph{stochastic} primal-dual (SPD) updates further reduces the computational effort by using one uncertainty realization per iteration to obtain:
\begin{subequations}
\begin{align}
     \bw_{k+1}&:=\big[\bw_{k}-\mu\sum_{n\in \mcN_g}\nabla_{\bw}c_n^k
     +\mu\blambda_k^\top\nabla_{\bw}\tilde{\mathbbm{1}}_{\epsilon}\left(\by^k-\bby\right)\big]_{\bpi_{\bw}\in\mcX}\label{eq:spd:p}\\
     \blambda_{k+1}&:=\big[\blambda_{k}+\nu(1-\alpha)\mathbf{1}-\nu\tilde{\mathbbm{1}}_{\epsilon}\left(\by^k-\bar{\by}\right)\big]_+,\label{eq:spd:dhigh}
\end{align}
\end{subequations}
where $c_n^k := c_n\left(\bpi_{\bw},\bphi^k\right)$ and $\by^k := \by\left(\bpi_{\bw},\bphi^k\right)$.

\subsection{Computing Gradients for Primal Updates}\label{subsec:grad}
The stochastic dual update simply requires evaluating the indicator approximation $\tilde{\mathbbm{1}}(x)$ for the constraint functions $\by^k$. The stochastic primal update of \eqref{eq:spd:p} however is more involved as it requires evaluating the gradients $\{\nabla_{\bw}c_n^k\}$ and the Jacobian matrix $\nabla_{\bw}\tilde{\mathbbm{1}}_{\epsilon}\left(\by^k-\bby\right)$. We next explain how these gradients can be computed via the chain rule and the inverse function theorem. 

We commence with gradient $\nabla_{\bw}c_n$ of the generation cost functions for all but the generator at the reference bus:
\begin{equation}\label{eq:gradcost}
\left(\nabla_{\bw}c_n\right)^\top
=\frac{\d c_n}{\d p_n^g}\cdot \left(\nabla_{\bx}p^g_n\right)^{\top}\cdot\nabla_{\bw}\bx,\quad\forall n \in \mcN_ g\setminus\{1\}.
\end{equation}
From the definition of $\bx$, the gradient $\nabla_{\bx}p^n_g$ is simply a canonical vector. The Jacobian matrix $\nabla_{\bw}\bx$ contains the partial derivatives of the DNN outputs with respect to the DNN inputs, and can be computed by most deep learning platforms such as TensorFlow using gradient \emph{back-propagation} in a computationally efficient manner. 

Computing $\nabla_{\bw}c_1$ for the generation cost at the reference bus is less direct, yet still computationally efficient. It is less direct because $p_1^g$ is not part of the setpoints, and hence, the DNN output $\bx$. Nonetheless, it does depend on all other generator setpoints and in turn $\bx$ indirectly, through the power flow equations. If we introduce the vector of voltages in polar coordinates $\bu:=[v_1~\dots~v_N~\theta_2~\dots~ \theta_N]^{\top}$ excluding the fixed angle $\theta_1=0$, the sought gradient can be computed as
\begin{equation}\label{eq:gradrefcost}
\left(\nabla_{\bw}c_1\right)^\top=
\frac{\d c_1}{\d p_1^g}\cdot
\left(\nabla_{\bu}p_1^g\right)^{\top}\cdot
\nabla_{\bx}\bu\cdot\nabla_{\bw}\bx
\end{equation}
The Jacobian matrix $\nabla_{\bx}\bu$ required in \eqref{eq:gradrefcost} can be found by solving a system of linear equations. Let vector $\bz$ collect the active and reactive power demands at all load buses. Then, using the power flow equations in \eqref{eq:pf_inj}, vectors $\bx$ and $\bz$ can be abstractly expressed as functions of $\bu$ as
\begin{align}\label{eq:syseq}
\begin{bmatrix}
    \bx\\
    \bz
\end{bmatrix}=
\begin{bmatrix}
    \bg(\bu)\\
    \bdell(\bu)
\end{bmatrix}.
\end{align}
Regarding the mapping $\bg(\bu)$, recall that $\bx:=\{\{v_n\}_{n\in\mcN_ g},\{p_n^g\}_{n\in\mcN_ g\setminus \{1\}}\}$. Then, the first $|\mcN_g|$ entries of $\bg(\bu)$ are simply $\be_n^\top\bu$, where $\be_n$ is the $n$-th canonical vector. The second $|\mcN_g|$ entries of $\bg(\bu)$ follow from the power flow equations plus any possible load $p_n^\ell$ at the corresponding bus. Since $p_1^g$ is a dependent variable, it has not been included in $\bx$ or $\bz$. Differentiating either sides of \eqref{eq:syseq} with respect to $\bx$ yields:
\begin{equation}\label{eq:diff}
     \begin{bmatrix}
        \bI\\
        \bzero
    \end{bmatrix}=
    \begin{bmatrix}
        \nabla_{\bu}\bg\\
        \nabla_{\bu}\bdell
    \end{bmatrix}\nabla_{\bx}\bu   
\end{equation}
where $\bI$ is an identity matrix of  dimensions $2|\mcN_ g|-1$ and $\bzero$ is a matrix of all zeroes of dimension $2|\mcN_\ell|\times (2|\mcN_ g|-1)$. From~\eqref{eq:diff}, we can now compute the Jacobian matrix $\nabla_{\bx}\bu$ as
\begin{align}\label{eq:gradxu}
    \nabla_{\bx}\bu=\begin{bmatrix}
        \nabla_{\bu}\bg(\bu)\\
        \nabla_{\bu}\bdell(\bu)
    \end{bmatrix}^{-1}
     \begin{bmatrix}
        \bI\\
        \bzero
    \end{bmatrix}
\end{align}
where the matrix to be inverted here is $(2N-1)\times (2N-1)$.

The Jacobian matrix $\nabla_{\bw}\tilde{\mathbbm{1}}_{\epsilon}\left(\by-\bar{\by}\right)$ appearing in the primal update of \eqref{eq:spd:p} can be found in a similar application of the chain rule and the inverse function theorem as
\begin{equation}\label{eq:gradindhigh}
\nabla_{\bw}\tilde{\mathbbm{1}}_{\epsilon}\left(\by-\bar{\by}\right)= \nabla_{\by}\tilde{\mathbbm{1}}_{\epsilon}\left(\by-\bar{\by}\right)\cdot 
\nabla_{\bu}\by\cdot 
\nabla_{\bx}\bu\cdot 
\nabla_{\bw}{\bx}.
\end{equation}
Matrix $\nabla_{\by}\tilde{\mathbbm{1}}_{\epsilon}\left(\by-\bar{\by}\right)$ is diagonal with diagonal entries provided by \eqref{eq:gradlog}. Matrix $\nabla_{\bu}\by$ contains the partial derivatives of all constraint functions of interest with respect to voltage magnitudes and angles. Such derivatives can be readily computed from \eqref{eq:pf_inj}--\eqref{eq:flow_VA}. 

\begin{figure}[t]
    \centering
    \includegraphics[scale=0.45]{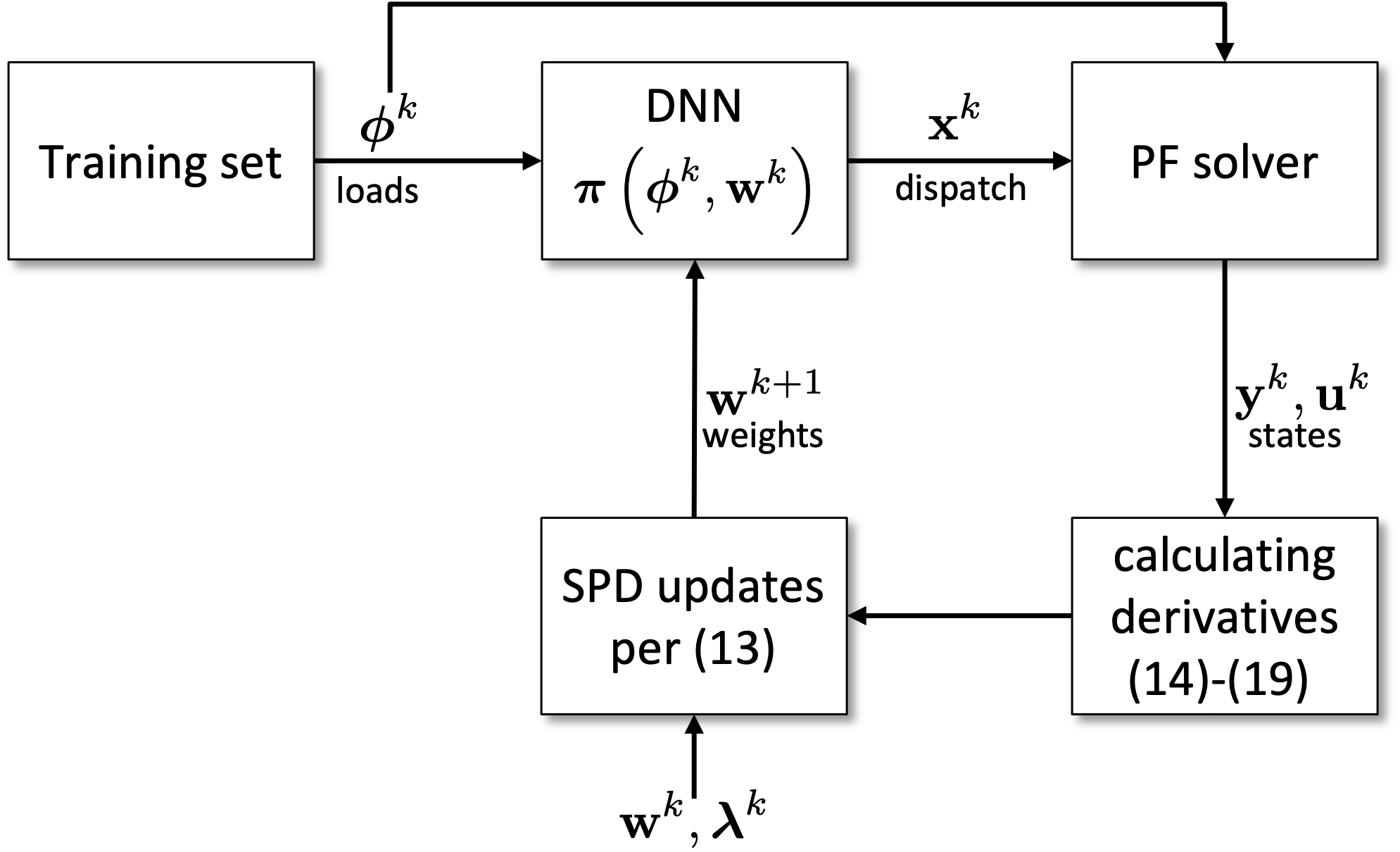}
    \caption{Overview of the training phase for the DNN-based policy.}
    \label{fig:training}
\end{figure}

\subsection{Deployment Workflow}\label{sec:overview}
Deploying the DNN-based strategy consists of two phases; (i) \emph{solving the SOPF}, which is achieved by the training phase of the DNN, and (ii) \emph{real time computation of control actions}, which is done by simply evaluating the trained DNN and is similar to the testing phase of DNNs. For the training phase, the system operator samples $K$ grid conditions $\{\bphi^k\}^{K}_{k=1}$ that reflect the real-time conditions to which the policy will be applied. The training samples could be sourced from historically recorded loads, simulations, or predictions. Depending on the generators that have been committed for the upcoming time-period, the operator identifies generator and load buses in sets $\mcN_ g$ and $\mcN_\ell$, respectively, and defines vectors $\bx$ and $\by$. The input and output layers of the DNN match the dimensions of $\bphi$ and $\bx$, accordingly. The dimensions and the number of hidden layers are hyper-parameters to be determined for the specific CC-OPF setting. The DNN is then trained using SPD as per Section~\ref{sec:DNNtrain}. This entails: \emph{i)} Sampling a $\bphi^k$; \emph{ii)} Performing a forward pass through the DNN to obtain $\bx^k$; \emph{iii)} Obtaining $\by^k$ and $\bu^k$ from a power flow solver; \emph{iv)} Using $\{\theta^k,\bx^k,\by^k,\bu^k\}$ to calculate the gradients \eqref{eq:gradcost}, \eqref{eq:gradrefcost}, and \eqref{eq:gradindhigh}; and \emph{v)} Performing updates \eqref{eq:spd:p}--\eqref{eq:spd:dhigh}. These steps are repeated over the training samples for multiple epochs. Figure~\ref{fig:training} depicts a block diagram of the training process. In real-time, the operator simply feeds the DNN with the real-time realizations of the uncertainty $\bphi$ to obtain the dispatch $\bx$ at the output.

\section{Numerical Tests} \label{sec:numerical}
\subsection{Experimental Setup}
The performance of the proposed DNN-based policy is evaluated using transmission networks of varying sizes. The simulation scripts are written in Python and run on a $2.4$~GHz 8-Core Intel Core i9 processor laptop computer with $64$~GB RAM. Compatibility with the commonly used MATPOWER~\cite{MATPOWER} models is achieved by interfacing the Python script with the open-source Octave engine~\cite{octave}. The communication between the two platforms is enabled via the Oct2py library~\cite{oct2py}, which allows seamless calling of M-files and Octave functions from Python. Leveraging upon this functionality, in-built MATPOWER functions are called to read the networks, solve the power flow equations, and calculate the derivatives needed for the primal updates. The  implementation advantages of such synergy come at the cost of the communication overhead between the two platforms. For the networks tested, this overhead is within $0.02-0.03$ seconds per iteration. Because these delays can be avoided by porting the required MATPOWER functions to Python, they have been eliminated from the training times reported here.

TensorFlow libraries are employed to model and train the DNN. For each of the experiments, a five-layer DNN is considered. The five layers are: \emph{1)} input layer matching the dimensions of $\bphi$; \emph{2)} two hidden layers, each with half the dimensions of the input layer; \emph{3)} output layer with the same dimensions as $\bx$; and \emph{4)} custom designed activation layer based on $\mathrm{tanh}$ ensuring $\bx\in\mcX$. All DNN weights are initialized by randomly drawing from a standard normal distribution. DNN biases and dual variables are initialized at zero. The primal updates are performed using the Adam optimizer and the dual variables are updated using SGD. The primal step size $\mu$ decays exponentially at the rate of $0.5$ per epoch, whereas the dual step size $\nu$ decays with the square-root of the iteration index~\cite{LKMG17}. The values for hyper-parameters $\{\epsilon,E,\mu_0,\nu_0\}$ are identified for each network and are reported with results from the experiment.

For training and testing the DNN, we generated $1,000$ samples of $\bphi$ by adding zero-mean uniformly distributed noise to the nominal loads of the MATPOWER models. Recall that training corresponds to solving the SOPF formulation in \eqref{eq:opf5} and testing corresponds to using the computed policy for real time prediction of control actions in response to uncertainty. To ensure that the problem in \eqref{eq:opf5} is feasible, the generated demands were truncated within a range of $[-R,R]$ around the nominal point, where $R$ was selected per network to ensure that the corresponding deterministic formulation of the OPF in \eqref{eq:opf2} is feasible. The $1,000$ generated samples were then randomly grouped into training and testing sets of $800$ and $200$ examples, respectively. Training was conducted over $E$ epochs of the training set. The performance of the DNN-based policy was bench-marked against the \emph{OPF-policy} that solves the OPF in \eqref{eq:opf2} deterministically per $\bphi_k$.

\subsection{Accuracy of the Logistic Function Approximation}
We investigate the accuracy and performance of the logistic function approximation using $\tilde{\mathbbm{1}}_{\epsilon}(x)$ to the chance constraints in \eqref{eq:opf4} using the $14$-bus ``pglib\_opf\_case14\_ieee'' network. For $\alpha\in\{0.05,0.10,0.15,0.20\}$ and $\{\epsilon,E,\mu_0,R\}=\{0.01,5,10^{-3},0.1\}$, Table~\ref{tab:14bus} shows the values of the remaining hyperparameters and the training time. The performance of the DNN-based policy on unseen test samples is also divulged via three metrics: maximum sampled probability of constraint violations (abbreviated in Table~\ref{tab:14bus} as \emph{maximum violation $[\%]$}), test time, and the average cost.

Table~\ref{tab:14bus} demonstrates the advantages of the DNN policy. First, the evaluated violation probabilities closely follows the prescribed $\alpha$ showing that the approximation $\tilde{\mathbbm{1}}_{\epsilon}(x)$ is sufficiently accurate while facilitating use of the efficient SPD algorithm. This is an improvement over the conservative convex approximation of the indicator function in \cite{OPFandLearnTSG21}, where the observed violation probabilities were considerably less than $\alpha$.
Second, during the real-time evaluation of the policy, the DNN-based policy is able to predict the dispatch for the $200$ test samples in around $0.3$ seconds. Comparing this to the OPF-policy, which has an evaluation time of $31.3$ seconds and cost of $\$2180.16$, the DNN policy is approximately $100$ times faster without sacrificing optimality. 

\begin{table}[t]
\centering
\caption{\label{tab:14bus} Training and testing details for the $14$-bus system for different values of $\alpha$.}\vspace{-.4em}
\resizebox{\columnwidth}{!}{%
\begin{tabular}{|r|r|r|r|r|r|r|}
\hline
\textbf{$\alpha$} & \multicolumn{1}{c|}{$\nu_0$} & \multicolumn{1}{c|}{\begin{tabular}[c]{@{}c@{}}Train\\ time (sec)\end{tabular}}&
\begin{tabular}[c]{@{}c@{}}Maximum \\violation $[\%]$ \end{tabular}& \begin{tabular}[c]{@{}c@{}} Test \\time (sec)\end{tabular} & \begin{tabular}[c]{@{}c@{}}Average\\cost $[\$]$\end{tabular} \\ \hline
$0.05$  &$3\cdot 10^{-4}$ & $97.4$&$3.5$ & $0.30$& $2180.53$\\ \hline
$0.10$ &$1.5\cdot 10^{-4}$ & $98.5$&$8.5$ & $0.33$ & $2180.48$ \\ \hline
$0.15$ &$1\cdot 10^{-4}$ & $100.6$& $16.5$ &  $0.34$ & $2180.45$\\ \hline
$0.20$ &$1.8\cdot 10^{-4}$ & $100.6$& $17.0$  & $0.34$ & $2180.44$ \\ \hline
\end{tabular}
}
\vspace{-2em}
\end{table}


\subsection{Enforcing Communication and Complexity Constraints}
While our presentation thus far has been focused on the unconstrained DNN policy where all controllable variables are allowed to respond to uncertainty realizations, the framework also allows for seamless encoding of communication constraints where only part of $\bphi$'s is observed and communicated, and complexity constraints where only a subset of the control variables respond to uncertainty. This flexibility allows us to account for any limitations of the available infrastructure.

As an example, we consider an automatic generation control (AGC) type policy, where only the active power generation of the participating generators responds to the total active load deviation in the system. This policy closely resembles the usual affine policy model in terms of input-output dependencies but allows for more general non-linear functions. This policy can be implemented within the DNN framework by modifying the input layer to have only $1$ neuron that receives the signal $\sum_{n \in \mcN}p^d_n$, and making the neurons for $\{v^g_n\}_{n \in \mcN^g}$ insensitive to the input by setting the corresponding weights to $0$. Note that the modified DNN still produces the average nominal values for $\{v^g_n\}_{n \in \mcN^g}$ at the output, because the biases in the output layer are learnt during training.

The AGC-type policy along with the full policy are evaluated on the $14$ bus system for $\alpha=0.1$ and two values of $R=\{0.1,0.2\}$ corresponding to a \emph{moderate} case with smaller uncertainty and a \emph{stressed} case with larger uncertainty respectively. The results are presented in Table~\ref{tab:test_agc}. For the moderate case, both policies are able to enforce the chance constraints with the full policy marginally out-performing the AGC-type policy in terms of cost. For the stressed case however, the AGC-type policy could not converge to anywhere near the desired $\alpha=0.10$ even after 20 epochs and saturated at a constraint violation probability of around $26.5\%$. The full policy on the other hand was able to decrease the probability of constraint violations to $10.5\%$. As a baseline comparison, the OPF policy was evaluated and found to attain an average cost of $\$~2183.14$. The stressed case demonstrates that general control policies can significantly improve cost of control actions and their ability to enforce constraints in real time.


\begin{table}[t]
\centering
\caption{\label{tab:test_agc} Test results comparing the full policy with an AGC-type policy for $\alpha=0.10$, and $E=20$.}
\vspace{-1em}
\begin{tabular}{cllll}
\multicolumn{1}{l}{}&  &    & &                                              \\ \hline
\multicolumn{1}{|c|}{\multirow{2}{*}{Policy}} & \multicolumn{2}{c|}{Maximum violation $[\%]$}  & \multicolumn{2}{c|}{Average cost $[\$]$} \\ \cline{2-5} 
\multicolumn{1}{|c|}{}  & \multicolumn{1}{l|}{\begin{tabular}[c]{@{}l@{}}$R=0.1$\end{tabular}} & \multicolumn{1}{l|}{\begin{tabular}[c]{@{}l@{}}$R=0.2$\end{tabular}} & \multicolumn{1}{l|}{\begin{tabular}[c]{@{}l@{}}$R=0.1$\end{tabular}} & \multicolumn{1}{l|}{\begin{tabular}[c]{@{}l@{}}$R=0.2$\end{tabular}} \\ \hline
\multicolumn{1}{|c|}{Full} & \multicolumn{1}{l|}{$8.5$}                   & \multicolumn{1}{l|}{$10.5$}  & \multicolumn{1}{l|}{$2180.45$}          & \multicolumn{1}{l|}{$2184.67$} \\ \hline
\multicolumn{1}{|c|}{AGC-type}& \multicolumn{1}{l|}{$9.0$} & \multicolumn{1}{l|}{$26.5$}  & \multicolumn{1}{l|}{$2185.52$}            & \multicolumn{1}{l|}{$2189.23$}   \\ \hline
\end{tabular}
\vspace{-1em}
\end{table}

\subsection{Scalability}
We test the DNN framework on a number of test networks of varying sizes.
The values for hyper-parameters, load variation parameter $R$, and the training times of these networks are compiled in Table~\ref{tab:train_nets} and Table~\ref{tab:test_nets} reports the performance of the trained DNNs. The testing times and the average costs attained by the proposed strategy are compared against the OPF policy. 
For all networks, the DNN policy is consistently able to enforce the chance constraints with the specified $\alpha$ with the costs remaining close to that of the OPF policy indicating that the DNN policy is near-optimal. The training times in Table~\ref{tab:train_nets} show a moderate increase with size of the test networks indicating that the proposed approach is highly scalable. We remark that the exact training times reported are less indicative than the trend. The exact times are highly dependent on implementation details where parallelized implementations and GPU deployment can drastically improve these numbers. 
The evaluation times in Table~\ref{tab:test_nets}  are orders of magnitude faster than the OPF policy, making them highly suitable for real-time computation of control actions.


\begin{table}[t]
\centering
\caption{\label{tab:train_nets} The values for hyper-parameters $\{\epsilon,E,\mu_0,\nu_0\}$, setting $R$, and the training times for different networks for $\alpha=0.05$ }
\vspace{-0.5em}
\begin{tabular}{|r|r|r|r|r|r|r|}
\hline
 \begin{tabular}[c]{@{}c@{}}Network\end{tabular}&\textbf{$\epsilon$} & \textbf{$E$} & \multicolumn{1}{c|}{$\mu_0$} & \multicolumn{1}{c|}{$\nu_0$} & \textbf{$R$} & \multicolumn{1}{c|}{\begin{tabular}[c]{@{}c@{}}Train \\ time (sec)\end{tabular}} \\ \hline
  case6ww&$0.005$ &$5$ & $10^{-3}$ &$4\cdot 10^{-2}$ & $0.05$ & $73.45$\\ \hline
    case69&$0.01$ &$5$ & $10^{-3}$ &$10^{-2}$ & $0.1$ & $79.85$\\ \hline
 case118&$0.07$ &$5$ & $10^{-3}$ &$1$ & $0.01$ & $308.4$\\ \hline
 case141&$0.01$ &$5$ & $10^{-3}$ &$10^{-3}$ & $0.1$ & $104.25$\\ \hline
\end{tabular}
\vspace{-0.5em}
\end{table}
\begin{table}[t]
\centering
\caption{\label{tab:test_nets}Test results for the full DNN-based policy (full) and OPF policy (OPF) for different networks for $\alpha=0.05$.}
\vspace{-1.5em}
\resizebox{\columnwidth}{!}{%
\begin{tabular}{cccllll}
\multicolumn{1}{l}{}& \multicolumn{1}{l}{}&    & & &              \\ \hline
\multicolumn{1}{|c|}{\multirow{2}{*}{\begin{tabular}[c]{@{}c@{}}Network\end{tabular}}} & \multicolumn{1}{c|}{\multirow{2}{*}{\begin{tabular}[c]{@{}c@{}}Maximum\\violation $[\%]$ \end{tabular}}}& \multicolumn{2}{c|}{Time (sec)}  & \multicolumn{2}{c|}{\begin{tabular}[c]{@{}c@{}}Average cost $[\$]$\end{tabular}} \\ \cline{3-6} 
\multicolumn{1}{|c|}{}  & \multicolumn{1}{c|}{} & \multicolumn{1}{l|}{Prop.} & \multicolumn{1}{l|}{Opt.} & \multicolumn{1}{l|}{Prop.} & \multicolumn{1}{l|}{Opt.} \\ \hline
\multicolumn{1}{|c|}{case6ww}&\multicolumn{1}{c|}{$5$} & \multicolumn{1}{c|}{ $0.22$} & \multicolumn{1}{c|}{$24.09$} &\multicolumn{1}{c|}{ $3.17\cdot 10^{3}$} &\multicolumn{1}{c|}{ $3.15\cdot 10^{3}$}\\ \hline
\multicolumn{1}{|c|}{case69}&\multicolumn{1}{c|}{$0$} & \multicolumn{1}{c|}{ $0.25$} & \multicolumn{1}{c|}{$26.35$} &\multicolumn{1}{c|}{ $80.58$} &\multicolumn{1}{c|}{ $80.58$}\\ \hline
\multicolumn{1}{|c|}{case118}&\multicolumn{1}{c|}{$0.5$} & \multicolumn{1}{c|}{ $0.39$} & \multicolumn{1}{c|}{$61.10$} &\multicolumn{1}{c|}{ $1.30\cdot 10^5$} &\multicolumn{1}{c|}{ $1.30\cdot 10^5$}\\ \hline
\multicolumn{1}{|c|}{case141}&\multicolumn{1}{c|}{$0.0$} & \multicolumn{1}{c|}{ $0.28$} & \multicolumn{1}{c|}{$40.34$} &\multicolumn{1}{c|}{$251.89$} &\multicolumn{1}{c|}{$251.89$}\\ \hline
\end{tabular}
}\vspace{-1em}
\end{table}

\section{Conclusions}\label{sec:conclusion}
We presented a DNN-based approach for solving the SOPF, where DNNs are used to parameterize the control policies required for real-time power balancing in response to uncertainty. Our approach does not require previously generated training labels and instead used the training phase to solve the SOPF. Stochastic primal-dual updates are employed to learn the DNN weights such that generation costs are minimized while respecting the power system constraints. Numerical tests on a variety of benchmark networks confirm that the generalized policy is able to provide high quality feasible  solutions to chance constrained AC-OPF problem over a range of operating conditions, with significant improvements over an AGC-type policy in terms of cost and constraint enforcement. Comparison with the OPF policy where an OPF is solved in response to each uncertainty realization shows that the DNN policy is able to attain similar levels of feasibility and optimality, while facilitating near-instant computation of real-time control actions. In future, we plan to extend our approach to joint chance-constrained OPF problems and research scenario selection policies to aid in faster DNN training. 

\bibliographystyle{IEEEtran}
\bibliography{myabrv,inverters,kekatos,power}

\begin{thebibliography}{10}
\providecommand{\url}[1]{#1}
\csname url@samestyle\endcsname
\providecommand{\newblock}{\relax}
\providecommand{\bibinfo}[2]{#2}
\providecommand{\BIBentrySTDinterwordspacing}{\spaceskip=0pt\relax}
\providecommand{\BIBentryALTinterwordstretchfactor}{4}
\providecommand{\BIBentryALTinterwordspacing}{\spaceskip=\fontdimen2\font plus
\BIBentryALTinterwordstretchfactor\fontdimen3\font minus
  \fontdimen4\font\relax}
\providecommand{\BIBforeignlanguage}[2]{{%
\expandafter\ifx\csname l@#1\endcsname\relax
\typeout{** WARNING: IEEEtran.bst: No hyphenation pattern has been}%
\typeout{** loaded for the language `#1'. Using the pattern for}%
\typeout{** the default language instead.}%
\else
\language=\csname l@#1\endcsname
\fi
#2}}
\providecommand{\BIBdecl}{\relax}
\BIBdecl

\bibitem{bienstock2014chance}
D.~Bienstock, M.~Chertkov, and S.~Harnett, ``Chance-constrained optimal power
  flow: Risk-aware network control under uncertainty,'' \emph{Siam Review},
  vol.~56, no.~3, pp. 461--495, 2014.

\bibitem{roald2013analytical}
L.~Roald, F.~Oldewurtel, T.~Krause, and G.~Andersson, ``Analytical
  reformulation of security constrained optimal power flow with probabilistic
  constraints,'' in \emph{2013 IEEE Grenoble Conference}.\hskip 1em plus 0.5em
  minus 0.4em\relax IEEE, 2013, pp. 1--6.

\bibitem{vrakopoulou2013probabilistic}
M.~Vrakopoulou, K.~Margellos, J.~Lygeros, and G.~Andersson, ``A probabilistic
  framework for reserve scheduling and security assessment of systems with high
  wind power penetration,'' \emph{{IEEE} Trans. Power Syst.}, vol.~28, no.~4,
  pp. 3885--3896, 2013.

\bibitem{dall2017chance}
E.~Dall’Anese, K.~Baker, and T.~Summers, ``Chance-constrained {AC} optimal
  power flow for distribution systems with renewables,'' \emph{{IEEE} Trans.
  Power Syst.}, vol.~32, no.~5, pp. 3427--3438, 2017.

\bibitem{MEZGHANI2020}
I.~Mezghani, S.~Misra, and D.~Deka, ``Stochastic {AC} optimal power flow: A
  data-driven approach,'' \emph{Electric Power Systems Research}, vol. 189, Dec
  2020.

\bibitem{muhlpfordt2019chance}
T.~M{\"u}hlpfordt, L.~Roald, V.~Hagenmeyer, T.~Faulwasser, and S.~Misra,
  ``Chance-constrained {AC} optimal power flow: A polynomial chaos approach,''
  \emph{{IEEE} Trans. Power Syst.}, vol.~34, no.~6, 2019.

\bibitem{metivier2020efficient}
D.~M{\'e}tivier, M.~Vuffray, and S.~Misra, ``Efficient polynomial chaos
  expansion for uncertainty quantification in power systems,'' \emph{Electric
  Power Systems Research}, vol. 189, p. 106791, 2020.

\bibitem{roald2015optimal}
L.~Roald, S.~Misra, M.~Chertkov, and G.~Andersson, ``Optimal power flow with
  weighted chance constraints and general policies for generation control,'' in
  \emph{Proc. {IEEE} Conf. on Decision and Control}.\hskip 1em plus 0.5em minus
  0.4em\relax IEEE, 2015, pp. 6927--6933.

\bibitem{lorca2017robust}
A.~Lorca and X.~Sun, ``{The Adaptive Robust Multi-Period Alternating Current
  Optimal Power Flow Problem},'' \emph{{IEEE} Trans. Power Syst.}, vol.~33,
  no.~2, pp. 1993--2003, 2018.

\bibitem{lee2021robust}
D.~Lee, K.~Turitsyn, D.~K. Molzahn, and L.~Roald, ``Robust {AC} optimal power
  flow with robust convex restriction,'' \emph{{IEEE} Trans. Power Syst.},
  2021.

\bibitem{JKGD19}
M.~Jalali, V.~Kekatos, N.~Gatsis, and D.~Deka, ``Designing reactive power
  control rules for smart inverters using support vector machines,''
  \emph{{IEEE} Trans. Smart Grid}, vol.~11, no.~2, pp. 1759--1770, Mar. 2020.

\bibitem{ng2018statistical}
Y.~Ng, S.~Misra, L.~A. Roald, and S.~Backhaus, ``Statistical learning for dc
  optimal power flow,'' in \emph{2018 Power Systems Computation Conference
  (PSCC)}.\hskip 1em plus 0.5em minus 0.4em\relax IEEE, 2018, pp. 1--7.

\bibitem{DekaMisraPowerTech19}
D.~{Deka} and S.~{Misra}, ``Learning for {DC-OPF}: {C}lassifying active sets
  using neural nets,'' in \emph{{IEEE PowerTech}}, Milan, Italy, Jun. 2019, pp.
  1--6.

\bibitem{zamzam2020learning}
A.~S. Zamzam and K.~Baker, ``Learning optimal solutions for extremely fast ac
  optimal power flow,'' in \emph{Proc. {IEEE} Intl. Conf. on Smart Grid
  Commun.}\hskip 1em plus 0.5em minus 0.4em\relax IEEE, 2020, pp. 1--6.

\bibitem{chen2020learning}
Y.~Chen and B.~Zhang, ``Learning to solve network flow problems via neural
  decoding,'' \emph{arXiv preprint arXiv:2002.04091}, 2020.

\bibitem{SGKCB2020}
M.~K. Singh, S.~Gupta, V.~Kekatos, G.~Cavraro, and A.~Bernstein, ``Learning to
  optimize power distribution grids using sensitivity-informed deep neural
  networks,'' in \emph{Proc. {IEEE} Intl. Conf. on Smart Grid Commun.}, Tempe,
  AZ, Nov. 2020, pp. 1--6.

\bibitem{L2OPWRS21}
\BIBentryALTinterwordspacing
M.~K. Singh, V.~Kekatos, and G.~B. Giannakis, ``Learning to solve the {AC-OPF}
  using sensitivity-informed deep neural networks,'' \emph{{IEEE} Trans. Power
  Syst.}, Jul. 2021, revised. [Online]. Available:
  \url{https://arxiv.org/abs/2103.14779}
\BIBentrySTDinterwordspacing

\bibitem{SG_SGC20}
S.~Gupta, V.~Kekatos, and M.~Jin, ``Deep learning for reactive power control of
  smart inverters under communication constraints,'' in \emph{Proc. {IEEE}
  Intl. Conf. on Smart Grid Commun.}, Tempe, AZ, 2020, pp. 1--6.

\bibitem{OPFandLearnTSG21}
\BIBentryALTinterwordspacing
------, ``Controlling smart inverters using proxies{: A} chance-constrained
  {DNN}-based approach,'' \emph{{IEEE} Trans. Smart Grid}, Jul. 2021,
  (revised). [Online]. Available: \url{https://arxiv.org/abs/2105.00429}
\BIBentrySTDinterwordspacing

\bibitem{ergodic}
G.~Wang, V.~Kekatos, A.-J. Conejo, and G.~B. Giannakis, ``Ergodic energy
  management leveraging resource variability in distribution grids,''
  \emph{{IEEE} Trans. Power Syst.}, vol.~31, no.~6, Nov. 2016.

\bibitem{6965430}
V.~Kekatos, G.~Wang, and G.~B. Giannakis, ``Stochastic loss minimization for
  power distribution networks,'' in \emph{Proc. North American Power
  Symposium}, Pullman, WA, Sep. 2014, pp. 1--6.

\bibitem{Ribeiro19}
M.~{Eisen}, C.~{Zhang}, L.~F.~O. {Chamon}, D.~D. {Lee}, and A.~{Ribeiro},
  ``Learning optimal resource allocations in wireless systems,'' \emph{{IEEE}
  Trans. Signal Processing}, vol.~67, no.~10, pp. 2775--2790, May 2019.

\bibitem{Nemirovski07}
A.~Nemirovski and A.~Shapiro, ``Convex approximations of chance constrained
  programs,'' \emph{SIAM Journal on Optimization}, vol.~17, no.~4, pp.
  969--996, 2007.

\bibitem{Chen95}
C.~Chen and O.~L. Mangasarian, ``Smoothing methods for convex inequalities and
  linear complementarity problems,'' \emph{Mathematical Programming}, vol.~71,
  no.~1, pp. 51--69, 1995.

\bibitem{MATPOWER}
R.~D. Zimmerman, C.~E. Murillo-Sanchez, and R.~J. Thomas, ``{MATPOWER}:
  steady-state operations, planning and analysis tools for power systems
  research and education,'' \emph{{IEEE} Trans. Power Syst.}, vol.~26, no.~1,
  pp. 12--19, Feb. 2011.

\bibitem{octave}
\BIBentryALTinterwordspacing
J.~W. Eaton, D.~Bateman, S.~Hauberg, and R.~Wehbring, ``{GNU Octave} version
  6.1.0 manual: a high-level interactive language for numerical computations,''
  2020. [Online]. Available:
  \url{https://www.gnu.org/software/octave/doc/v6.1.0/}
\BIBentrySTDinterwordspacing

\bibitem{oct2py}
\BIBentryALTinterwordspacing
``Oct2py: Python to gnu octave bridge.'' [Online]. Available:
  \url{https://oct2py.readthedocs.io/en/latest/index.html}
\BIBentrySTDinterwordspacing

\bibitem{LKMG17}
L.~M. Lopez-Ramos, V.~Kekatos, A.~G. Marques, and G.~B. Giannakis,
  ``Two-timescale stochastic dispatch of smart distribution grids,''
  \emph{{IEEE} Trans. Smart Grid}, vol.~9, no.~5, pp. 4282--4292, Sep. 2018.

\end{thebibliography}
\end{document}